\documentclass[final,3p]{elsarticle}

\usepackage{mathptmx}

\usepackage[numbers]{natbib}

\biboptions{sort&compress}

\usepackage{hyperref}

\usepackage{amsmath}
\usepackage{amsfonts}
\usepackage{amsthm}
\usepackage{mathtools}

\journal{Applied Numerical Mathematics}

\numberwithin{equation}{section}

\begin{document}

\begin{frontmatter}


  \title{Summation of divergent power series by means of factorial
    series}

\author{Ernst Joachim Weniger}
\ead{joachim.weniger@chemie.uni-regensburg.de}
\address{Institut f\"ur Physikalische und Theoretische Chemie,
  Universit\"at Regensburg, D-93040 Regensburg, Germany}

\begin{abstract}
  Factorial series played a major role in Stirling's classic book
  \emph{Methodus Differentialis} (1730), but now only a few specialists
  still use them. This article wants to show that this neglect is
  unjustified, and that factorial series are useful numerical tools for
  the summation of divergent (inverse) power series.  This is documented
  by summing the divergent asymptotic expansion for the exponential
  integral $E_{1} (z)$ and the factorially divergent
  Rayleigh-Schr\"{o}dinger perturbation expansion for the quartic
  anharmonic oscillator. Stirling numbers play a key role since they
  occur as coefficients in expansions of an inverse power in terms of
  inverse Pochhammer symbols and vice versa. It is shown that the
  relationships involving Stirling numbers are special cases of more
  general orthogonal and triangular transformations.
\end{abstract}

\begin{keyword} Factorial series; Divergent asymptotic (inverse) power
  series; Stieltjes series; Quartic anharmonic oscillator; Stirling
  numbers; General orthogonal and triangular transformations;

\MSC[2010] 11B73; 40A05; 40G99; 81Q15;

\end{keyword}

\end{frontmatter}


\typeout{==> Section: Introduction}
\section{Introduction}
\label{Sec:Introduction}

Power series are extremely important analytical tools not only in
mathematics, but also in the mathematical treatment of scientific and
engineering problems. Unfortunately, a power series representation for a
given function is from a numerical point of view a mixed blessing. A
power series converges within its circle of convergence and diverges
outside. Circles of convergence normally have finite radii, but there are
many series expansions of considerable practical relevance, for example
asymptotic expansions for special functions or quantum mechanical
perturbation expansions, whose circles of convergence shrink to a single
point.

The summation of divergent (inverse) power series is an old problem of
mathematics, which is of considerable relevance also in related
disciplines.  Many different summation techniques have been developed
which are often capable of associating a finite value to a divergent
series. A highly condensed overview of various summation techniques was
recently given in \cite[Appendices A and B]{Weniger/2008}. The role of
summation techniques in theoretical physics was discussed in the recent
review \cite{Caliceti/Meyer-Hermann/Ribeca/Surzhykov/Jentschura/2007}.

The topic of this article is the summation of divergent power series via
so-called \emph{factorial series}. A factorial series for a function
$\Omega \colon \mathbb{C} \to \mathbb{C}$, which vanishes as $z \to +
\infty$, is an expansion of the following type:
\begin{equation}
  \label{DefFactSer}
\Omega (z) \; = \; \frac {a_0} {z} \, + \, \frac {a_1 1!} {z (z+1)}
\, + \, \frac {a_2 2!} {z (z+1) (z+2)} \, + \, \cdots \; = \;
\sum_{\nu=0}^{\infty} \frac {a_{\nu} {\nu}!} {(z)_{\nu+1}} \, .
\end{equation}
Here, $(z)_{\nu+1} = \Gamma (z+\nu+1) / \Gamma (z) = z (z+1) \ldots
(z+\nu)$ is a Pochhammer symbol. The function $\Omega (z)$ represented by
the factorial series (\ref{DefFactSer}) may have simple poles at $z = -
m$ with $m \in \mathbb{N}_0$. The definition (\ref{DefFactSer}) is
typical of the mathematical literature. It will become clear later that
the separation of the series coefficients into a factorial $n!$ and a
reduced coefficient $a_n$ often offers formal advantages.

The use of factorial series for the summation of divergent (inverse)
power series is not a new idea. There is Watson's classic article on the
transformation of an asymptotic inverse power series to a convergent
factorial series \cite{Watson/1912b}. This topic was also considered in
articles by Nevanlinna \cite{Nevanlinna/1919}, Malgrange
\cite{Malgrange/1995}. and Ramis and Thomann \cite{Ramis/Thomann/1981}.
Thomann \cite{Thomann/1990,Thomann/1995} discussed the summation of
formal power series with the help of generalizations of factorial series,
where the Pochhammer symbols $(z)_{n+1} = z(z+1) \dots (z+n)$ are
replaced by products $z (z+t_{1}) \dots (z+t_{n})$. There is also a
recent article by Delabaere and Rasoamanana
\cite{Delabaere/Rasoamanana/2007} on the connection of Borel summation
and factorial series, which can be viewed to be an extension of a short
discussion in Borel's classic book \cite[pp.\ 234 - 245]{Borel/1988}.

These references show that there is no lack of knowledge about factorial
series in general and about their use as summation tools. Unfortunately,
this knowledge is restricted to a relatively small group of specialists,
and there is a deplorable lack of public awareness about factorial series
and their numerical usefulness. My claim is supported by the fact that
factorial series are not mentioned in the 2010 Mathematics Subject
Classification (MSC), in the web database MathWorld
({\verb+http://mathworld.wolfram.com+}), or in Wikipedia
({\verb+http://en.wikipedia.org/wiki/Main_Page+}). This neglect is not
justified, and I believe that the full potential of factorial series as
summation tools has not yet been realized.

Section \ref{Sec:RediscoveryFactSer} describes how I had become
interested in factorial series. In Section
\ref{Sec:BasicPropertiesFactSer}, the basic properties of factorial
series are reviewed. In Section \ref{Sec:InvPowSer<->FactSer}, it is
shown that inverse power series and factorial series can be transformed
into each other with the help of Stirling numbers. In Section
\ref{Sec:StieltjesSeries}, the transformation of Stieltjes series is
considered. Section \ref{Sec:SumAsySer_ExpInt} describes the
transformation of the divergent asymptotic series for the exponential
integral to a convergent factorial series. In Section
\ref{Sec:PowSer2FactSer}, a transformation for power series, which is
essentially a factorial series in $1/z$, is discussed, and in Section
\ref{Sec:TheQuarticAnharmonicOscillator}, this transformation formula is
used for the summation of the divergent Rayleigh-Schr\"{o}dinger
perturbation expansions for the ground state energy eigenvalue of the
quartic anharmonic oscillator. This article is concluded by a short
outlook in Section \ref{Sec:Outlook}. Those properties of Stirling
numbers, which are for our purposes most relevant, are reviewed in
\ref{App:StirlingNumbers}. In
\ref{App:SequenceInversionByTriangularMatrices}, it is shown that the
transformation formulas for inverse power and factorial series considered
in Section \ref{Sec:StieltjesSeries} are just special cases of more
general transformation formulas involving triangular and orthogonal
matrices.

\typeout{==> Section: My ``rediscovery'' of factorial series}
\section{My ``rediscovery'' of factorial series}
\label{Sec:RediscoveryFactSer}

In 1985/1986 I became interested in Levin's sequence transformation
\cite{Levin/1973}, and I tried to understand the mathematical theory
behind it. My interest was aroused by two articles by Smith and Ford
\cite{Smith/Ford/1979,Smith/Ford/1982} who had shown that certain
variants of Levin's sequence transformation were among the most powerful
as well as most versatile sequence transformations known at that time.

In my work on Levin's sequence transformation \cite{Levin/1973} I
discovered that its derivation becomes almost trivially simple if we
start from the model sequence \cite[Eq.\ (3.2-9)]{Weniger/1989}
\begin{equation}
  \label{Mod_Seq_Om}
  s_n \; = \; s \, + \, \omega_n z_n \, , \qquad n \in \mathbb{N}_0 \, .
\end{equation}
The remainder estimates $\omega_n$ are assumed to be known, and the
correction terms $z_n$ should be chosen in such a way that the products
$\omega_n z_n$ provide sufficiently accurate and rapidly convergent
approximations to the actual remainders $\{ r_n \}_{n=0}^{\infty}$ of the
sequence $\{ s_n \}_{n=0}^{\infty}$ which is to be transformed. 

In this approach, only the correction terms $\{ z_n \}_{n=0}^{\infty}$
have to be determined. If good remainder estimates can be found, the
determination of $z_n$ and the subsequent elimination of $\omega_n z_n$
from $s_n$ often leads to better results than the construction and
subsequent elimination of other approximations to $r_n$.

The model sequence (\ref{Mod_Seq_Om}) has another indisputable advantage:
A sequence transformation, which is exact for this model sequence, can be
constructed easily under very mild conditions. Let us assume that a
\emph{linear} operator $\hat{T}$ can be found which annihilates the
correction term $z_n$ for all $n \in \mathbb{N}_0$ according to $\hat{T}
(z_n) = 0$. Then we obtain a sequence transformation, which is exact for
the model sequence (\ref{Mod_Seq_Om}), by applying $\hat{T}$ to the ratio
$[s_n - s] / \omega_n = z_n$. Since $\hat{T}$ annihilates $z_n$ and is
by assumption linear, the following sequence transformation $\mathcal{T}$
is exact for the model sequence (\ref{Mod_Seq_Om}) \cite[Eq.\
(3.2-11)]{Weniger/1989}:
\begin{equation}
  \label{GenSeqTr}
  \mathcal{T} (s_n, \omega_n) \; = \; \frac
  {\hat{T} (s_n / \omega_n )} {\hat{T} (1 / \omega_n )} \; = \; s \, .
\end{equation}

The annihilation operator approach was introduced in \cite[Section
3.2]{Weniger/1989} in connection with my rederivation of Levin's
transformation \cite{Levin/1973}, but I also constructed in this way some
other, closely related sequence transformations \cite[Sections 7 -
9]{Weniger/1989}. Brezinski and Redivo Zaglia
\cite{Brezinski/RedivoZaglia/1994a,Brezinski/RedivoZaglia/1994b} and
Brezinski and Matos \cite{Brezinski/Matos/1996} showed later that this
approach is actually much more general: The majority of the currently
known sequence transformations can be derived via (\ref{GenSeqTr}) (for
further references on this topic, see \cite[p.\ 1214]{Weniger/2004}).

We obtain a model sequence for Levin's sequence transformation
\cite{Levin/1973} by assuming that $z_{n}$ in (\ref{Mod_Seq_Om}) is a
truncated power series in $1/(\beta+n)$ \cite[Eq.\
(7.1-1)]{Weniger/1989}:
\begin{equation}
  \label{z_n_Lev}
  z_{n} \; = \; \sum_{j=0}^{k-1} \, \frac{c_{j}}{(\beta+n)^{j}} \, ,
  \qquad \beta > 0 \, .
\end{equation}
The $k$th power of the finite difference operator $\Delta$ defined by
$\Delta f (n) = f (n+1) - f (n)$ annihilates an arbitrary polynomials
$P_{k-1} (n)$ of degree $k-1$ in $n$. Thus, the weighted difference
operator $\hat{T} = \Delta^k (\beta+n)^{k-1}$ is the appropriate
annihilation operator for $z_{n}$ defined by (\ref{z_n_Lev}), and Levin's
sequence transformation \cite{Levin/1973} can in the notation of
\cite[Eqs.\ (7.1-6) and (7.1-7)]{Weniger/1989} be expressed as follows:
\begin{equation}
  \label{GenLevTr}
  \mathcal{L}_{k}^{(n)} (\beta, s_n, \omega_n) \; = \; \frac
  {\Delta^k [(\beta+n)^{k-1} s_n / \omega_n]}
  {\Delta^k [(\beta+n)^{k-1} / \omega_n]} \; = \; \frac
  {\displaystyle
    \sum_{j=0}^{k} \, (-1)^{j} \, {\binom{k}{j}} \,
    \frac
    {(\beta + n +j )^{k-1}} {(\beta + n + k )^{k-1}} \,
    \frac {s_{n+j}} {\omega_{n+j}} }
  {\displaystyle
    \sum_{j=0}^{k} \, (-1)^{j} \, {\binom{k}{j}} \,
    \frac
    {(\beta + n +j )^{k-1}} {(\beta + n + k )^{k-1}} \,
    \frac {1} {\omega_{n+j}} }
  \, , \qquad k, n \in \mathbb{N}_0 \, .
\end{equation}
The numerator and denominator sums of this and of related transformations
can also be computed recursively \cite[Section III]{Weniger/2004}.

This undeniable success inspired me to look for other applications of the
annihilation operator approach. If we replace in (\ref{z_n_Lev}) the
powers $(\beta+n)^{j}$ by Pochhammer symbols $(\beta+n)_{j}$, we obtain a
truncated factorial series in $\beta+n$:
\begin{equation}
  \label{z_n_Wen}
  z_{n} \; = \; \sum_{j=0}^{k-1} \, \frac{c'_{j}}{(\beta+n)_{j}} \, ,
  \qquad \beta > 0 \, .
\end{equation}
Now, $\hat{T} = \Delta^k (\beta+n)_{k-1}$ is the appropriate annihilation
operator, and we obtain \cite[Eqs.\ (8.2-6) and (8.2-7)]{Weniger/1989}:
\begin{equation}
  \label{GenWenTr}
  \mathcal{S}_{k}^{(n)} (\beta, s_n, \omega_n) \; = \; \frac
  {\Delta^k [(\beta+n)_{k-1} s_n / \omega_n]}
  {\Delta^k [(\beta+n)_{k-1} / \omega_n]} \; = \; \frac
  {\displaystyle
    \sum_{j=0}^{k} \, (-1)^{j} \, {\binom{k}{j}} \,
    \frac
    {(\beta + n +j )_{k-1}} {(\beta + n + k )_{k-1}} \,
    \frac {s_{n+j}} {\omega_{n+j}} }
  {\displaystyle
    \sum_{j=0}^{k} \, (-1)^{j} \, {\binom{k}{j}} \,
    \frac
    {(\beta + n +j )_{k-1}} {(\beta + n + k )_{k-1}} \,
    \frac {1} {\omega_{n+j}} }
  \, , \qquad k, n \in \mathbb{N}_0 \, .
\end{equation}

My derivation of this sequence transformations, whose theory was
developed in \cite[Section 8]{Weniger/1989}, was entirely based on
heuristics. I had only looked for situations in which I could apply the
annihilation operator formalism (\ref{GenSeqTr}) effectively. If the
correction term $z_{n}$ is according to (\ref{z_n_Wen}) a truncated
factorial series, it can be annihilated easily, but I had no idea whether
the resulting sequence transformation $\mathcal{S}_{k}^{(n)} (\beta, s_n,
\omega_n)$ would be computationally useful or not.

$\mathcal{S}_{k}^{(n)} (\beta, s_n, \omega_n)$ was first used for the
evaluation of auxiliary functions in molecular electronic structure
calculations \cite{Weniger/Steinborn/1989a}. Later, it was used with
considerable success in the case of slowly convergent or divergent
alternating series (numerous references are listed in \cite[p.\
1225]{Weniger/2004}). Currently, $\mathcal{S}_{k}^{(n)} (\beta, s_n,
\omega_n)$ is used quite a lot in optics
\cite{Borghi/2007,Borghi/2008a,Borghi/2008b,Borghi/2008c,Borghi/2009,%
  Borghi/Alonso/2007,Borghi/Santarsiero/2003,%
  Li/Zang/Tian/2009,Li/Zang/Li/Tian/2009}.

When I constructed $\mathcal{S}_{k}^{(n)} (\beta, s_n, \omega_n)$ in
1986, I had no idea what the correction term (\ref{z_n_Wen}) actually is:
Factorial series had not been part of my mathematical training. It took a
while until I found out that (\ref{z_n_Wen}) is a truncated factorial
series, and that Nielsen's classic book \cite{Nielsen/1965}, which is
still one of the principal references on factorial series, had already
been waiting for quite a while on my bookshelf.

In my later work on convergence acceleration and summation processes, I
noticed that $\mathcal{L}_{k}^{(n)} (\beta, s_n, \omega_n)$ and
$\mathcal{S}_{k}^{(n)} (\beta, s_n, \omega_n)$ usually have similar, but
not identical properties. Nevertheless, in some cases spectacular
differences were observed. For example, in summation calculations for the
divergent Rayleigh-Schr\"{o}dinger perturbation expansions of the ground
state energies of anharmonic oscillators
\cite{Weniger/Cizek/Vinette/1993}, we observed that Levin's
transformation $\mathcal{L}_{k}^{(n)} (\beta, s_n, \omega_n)$ ultimately
produces divergent results, whereas $\mathcal{S}_{k}^{(n)} (\beta, s_n,
\omega_n)$ produces very good results (see also \cite{Weniger/1992}). A
similar divergence of Levin's transformation was also observed by
{\v{C}\'{\i}\v{z}ek}, Zamastil, and {Sk\'{a}la}
\cite{Cizek/Zamastil/Skala/2003}. Needless to say that these observations
puzzled me.

The construction of Levin's transformation is based on the assumption
that the ratio $[s_{n}-s]/\omega_{n}$ can be expressed as an inverse
power series, whereas $\mathcal{S}_{k}^{(n)} (\beta, s_n, \omega_n)$
implicitly assumes that $[s_{n}-s]/\omega_{n}$ can be expressed as a
factorial series. For me, it was a plausible hypothesis that the observed
differences of $\mathcal{L}_{k}^{(n)} (\beta, s_n, \omega_n)$ and
$\mathcal{S}_{k}^{(n)} (\beta, s_n, \omega_n)$ could be related to
different properties of inverse power and factorial series, respectively.
Because of my complete lack of knowledge about factorial series, I first
had to study their properties. The results presented in this article are
in some sense a by-product of these studies.

\typeout{==> Section: Basic properties of factorial series}
\section{Basic properties of factorial series}
\label{Sec:BasicPropertiesFactSer}

Factorial series have been known for a very long time. In Tweedle's
annotated translations of Stirling's classic \emph{Methodus
  Differentialis} it is remarked that Stirling was not the inventor of
factorial series. Apparently, Stirling became aware of factorial series
by the work of the French mathematician Nicole \cite[p.\
174]{Tweddle/2003}.  However, Stirling used factorial series extensively
and thus did a lot to popularize them.

The application of higher powers of the finite difference operator
$\Delta = \Delta_{z}$ to a factorial series in $z$ yields an extremely
compact result. If we use $\Delta_{z}^{k} [n!/(z)_{n+1}] = (-1)^k
(n+k)!/(z)_{n+k+1}$ with $k \in \mathbb{N}_0$, we obtain
\begin{equation}
  \label{Diff_k_FactSer}
  \Delta_{z}^{k} \, \Omega (z) \; = \; \sum_{\nu=0}^{\infty} \, 
  \Delta_{z}^{k} \, \frac {a_{\nu} {\nu}!} {(z)_{\nu+1}} \; = \; 
  (-1)^k \, \sum_{\nu=0}^{\infty} \, \frac {a_{\nu} (\nu+k)!}
  {(z)_{\nu+k+1}} \; = \; (-1)^k \, \sum_{\kappa=k}^{\infty} \, \frac
  {a_{\kappa-k} \kappa!} {(z)_{\kappa+1}} \, .  
\end{equation}
Factorial series play a similar role in the theory of difference
equations as inverse power series in the theory of differential
equations, and classic books on finite difference such as the ones by
Milne-Thomson \cite{Milne-Thomson/1981}, or N\"orlund
\cite{Noerlund/1926,Noerlund/1929,Noerlund/1954}) treat factorial series.
A contemporary discussion of the use of factorial series in the context
of difference equations can be found in a recent article by Olde Daalhuis
\cite{OldeDaalhuis/2004b}.

But I am are much more interested in the convergence properties of
factorial series, which fortunately can be analyzed easily. If we use
\cite[Eq.\ (6.1.47) on p.\ 257]{Abramowitz/Stegun/1972}
\begin{equation}
  \label{AsyGammaRatio}
\Gamma(z+a)/\Gamma(z+b) \; = \; z^{a-b} \,
\bigl[ 1 + \mathrm{O} (1/z) \bigr] \, , \qquad z \to \infty \, ,
\end{equation}
we obtain the asymptotic estimate $n!/(z)_{n+1} = \mathrm{O} (n^{-z})$ as
$n \to \infty$. Thus, the factorial series (\ref{DefFactSer}) converges
with the possible exception of the points $z = - m$ with $m \in
\mathbb{N}_0$ if and only if the associated Dirichlet series
$\tilde{\Omega} (z) = \sum_{n=1}^{\infty} a_n/n^z$ converges (see for
example \cite[p.\ 262]{Knopp/1964} or \cite[p.\ 167]{Landau/1906}).
Accordingly, a factorial series converges for sufficiently large $\Re
(z)$ even if the reduced series coefficients $a_{n}$ in
(\ref{DefFactSer}) grow like a fixed power $n^{\alpha}$ with $\alpha > 0$
as $n \to \infty$.

Factorial series are closely related to the beta function, which is
usually defined as the following ratio of gamma functions \cite[Eq.\
(6.2.2)]{Abramowitz/Stegun/1972}:
\begin{equation}
  \label{Def_BetaFun}
B (x, y) \; = \; \frac{\Gamma(x) \Gamma(y)}{\Gamma(x+y)} \, ,
\qquad x, y \in \mathbb{C} \, .
\end{equation}
Thus, the ratio $n!/(z)_{n+1}$ can be expressed as a beta function:
\begin{equation}
  \label{Beta2FSterm}
B (z, n+1) \; = \; \frac{n!}{(z)_{n+1}} \, .
\end{equation}
Accordingly, a factorial series can also be expressed as an expansion in
terms of beta functions \cite[p.\ 288]{Milne-Thomson/1981}:
\begin{equation}
  \label{FactSerBetaExpan}
\Omega (z) \; = \; \sum_{n=0}^{\infty} \, a_n \, B (z, n+1) \, .
\end{equation}

The beta function possesses numerous integral representations. For our
purposes the most useful one is the so-called \emph{Euler integral of the
  first kind} (see for example \cite[Eq.\
(6.2.1)]{Abramowitz/Stegun/1972}):
\begin{equation}
  \label{BetaFunIntRep}
B (x, y) \; = \; \int_{0}^{1} \, t^{x-1} \, (1-t)^{y-1} \, \mathrm{d} t
\, , \qquad \Re (x), \Re (y) > 0 \, .
\end{equation}

Combination of (\ref{FactSerBetaExpan}) and (\ref{BetaFunIntRep}) yields the
following integral representation:
\begin{equation}
  \label{IntRepFSterm}
\frac{n!}{(z)_{n+1}} \; = \; 
\int_{0}^{1} \, t^{z-1} \, (1-t)^{n} \, \mathrm{d} t
\, , \qquad \Re (z) > 0 \, \quad n \in \mathbb{N}_0 \, .
\end{equation}
If we now combine (\ref{Beta2FSterm}), (\ref{FactSerBetaExpan}) and
(\ref{IntRepFSterm}) and interchange integration and summation, we obtain
the following integral representation (see for instance \cite[Satz I on
p.\ 244]{Nielsen/1965} or \cite[p.\ 289]{Milne-Thomson/1981}):
\begin{subequations}
  \label{FS_IntRep}
\begin{align}
  \label{FS_IntRep_a}
  \Omega (z) & \; = \; \int_{0}^{1} \, t^{z-1} \, \varphi (t) \,
  \mathrm{d} t \, , \qquad \Re (z) > 0 \, ,
  \\
  \label{FS_IntRep_b}
  \varphi (t) & \; = \; \sum_{n=0}^{\infty} \, a_n \, (1-t)^n \, .
\end{align}
\end{subequations}
This integral representation is of considerable importance: Frequently,
the properties of $\Omega (z)$ can be studied more conveniently via the
corresponding properties of the conjugated function $\varphi (t)$ than
via the defining factorial series (\ref{DefFactSer}) (see for example
\cite[Kapitel XVII]{Nielsen/1965}). As discussed in Section
\ref{Sec:TheQuarticAnharmonicOscillator}, this integral representation
can also be used for the evaluation of factorial series.

\typeout{==> Section: Transformations of inverse power series and
  factorial series}
\section{Transformations of inverse power series and factorial series}
\label{Sec:InvPowSer<->FactSer}

Inverse powers $1/z^{k+1}$ and inverse Pochhammer symbols $1/(z)_{k+1}$
can be transformed into each other via (\ref{Poch_expand_Pow}) and
(\ref{InvPow_expand_InvPoch}), respectively. Therefore, inverse power
series and factorial series can also be transformed into each other. This
has been known for a very long time. The algebraic processes effecting
these transformations were already described in Nielsen's book
\cite{Nielsen/1965} which was first published in 1906. It seems, however,
that these potentially very useful transformation formulas are now
largely forgotten.

Let us assume that a function $\psi \colon \mathbb{C} \to \mathbb{C}$
possesses the inverse power series $\psi (z) = \sum_{n=0}^{\infty}
c_{n}/z^{n+1}$. If we insert (\ref{InvPow_expand_InvPoch}) into this
series and rearrange the order of summations, we obtain the
transformation formula
\begin{equation}
  \label{InvPowSer->FactSer}
\sum_{n=0}^{\infty} \, \frac{c_n}{z^{n+1}} \; = \;
\sum_{m=0}^{\infty} \, \frac{(-1)^{m}}{(z)_{m+1}} \, 
\sum_{\mu=0}^{m} \, (-1)^{\mu} \, \textbf{S}^{(1)} (m, \mu) \, 
c_{\mu} \, ,
\end{equation}
which shows that the coefficients of the factorial series are weighted
averages of the power series coefficients involving Stirling numbers of
the first kind. I employed this transformation already in
\cite{Weniger/2007a} in order to speed up the convergence of truncated
asymptotic expansions for the truncation errors of series expansions for
special functions.

An inverse expansion can also be derived. Let us assume that a function
$\chi \colon \mathbb{C} \to \mathbb{C}$ possesses a factorial series
$\chi (z) = \sum_{n=0}^{\infty} d_{n}/(z)_{n+1}$. If we insert
(\ref{Poch_expand_Pow}) into this series and rearrange the order of
summations, we obtain the transformation formula
\begin{equation}
  \label{FactSer->InvPowSer}
  \sum_{n=0}^{\infty} \, \frac{d_n}{(z)_{n+1}} \; = \;
  \sum_{m=0}^{\infty} \, \frac{(-1)^{m}}{z^{m+1}} \, 
  \sum_{\mu=0}^{m} \, (-1)^{\mu} \, 
  \textbf{S}^{(2)} (m, \mu) \, d_{\mu} \, .
\end{equation}

The operations producing the transformation formulas
(\ref{InvPowSer->FactSer}) and (\ref{FactSer->InvPowSer}) are purely
formal. Therefore, we cannot tacitly assume that the inverse power or the
factorial series necessarily converge. This has to be checked explicitly
in each case.

The sign patterns in the inner sums $\sum_{\mu=0}^{m} (-1)^{\mu}
\textbf{S}^{(1)} (m, \mu) c_{\mu}$ and $\sum_{\mu=0}^{m} (-1)^{\mu}
\textbf{S}^{(2)} (m, \mu) d_{\mu}$ in (\ref{InvPowSer->FactSer}) and
(\ref{FactSer->InvPowSer}) are of crucial importance for the convergence
or divergence of the formal expansions (\ref{InvPowSer->FactSer}) and
(\ref{FactSer->InvPowSer}). If the signs of the terms in the inner sums
alternate, we can hope for a substantial cancellation as in binomial sums
$\sum_{j=0}^{k} (-1)^{k} \binom{k}{j} f_{n+j}$, but if all terms have the
same sign, a potentially explosive accumulation can take place.  In the
former case, convergence is likely, while in the latter case we should be
prepared for divergence.

The ability of the Stirling numbers of the first kind to achieve a
cancellation is immediately obvious from its finite generating function.
Setting $z=k$ with $1 \le k \le n-1$ in (\ref{St1_GenFun1}) yields
$\sum_{\nu=0}^{n} k^{\nu} \textbf{S}^{(1)} (n, \nu) = 0$ for $n \ge 2$.
Let us now assume that the coefficients $c_{n}$ in
(\ref{InvPowSer->FactSer}) have strictly alternating signs. Consequently,
(\ref{St1_SignConv}) implies that the terms of the inner sum $(-1)^{m}
\sum_{\mu=0}^{m} (-1)^{\mu} \textbf{S}^{(1)} (m, \mu) c_{\mu} =
\sum_{\mu=0}^{m} \vert \textbf{S}^{(1)} (m, \mu) \vert c_{\mu}$ in
(\ref{InvPowSer->FactSer}) also have strictly alternating signs, and we
can hope for a substantial cancellation. The inner sum should for large
values of the outer index $m$ be much smaller in magnitude than its
individual terms. Because of this cancellation, we can hope that a
function defined by a \emph{divergent} inverse power series with strictly
alternating coefficients can be expressed and computed by a
\emph{convergent} factorial series.

The situation is not nearly as nice if the coefficients $c_{n}$ in
(\ref{InvPowSer->FactSer}) all have the same sign. If we set in
(\ref{St1_GenFun2}) $z=1$ and use (\ref{St1_SignConv}), we obtain
$(-1)^{m} \sum_{\mu=0}^{m} (-1)^{\mu} \textbf{S}^{(1)} (m, \mu) =
\sum_{\mu=0}^{m} \bigl\vert \textbf{S}^{(1)} (m, \mu) \bigr\vert = m!$.
Accordingly, the coefficients $\sum_{\mu=0}^{m} \bigl\vert
\textbf{S}^{(1)} (m, \mu) \bigr\vert c_{\mu}$ of the factorial series
should be (much) larger than the coefficients $c_{n}$ of the inverse
power series.

\typeout{==> Section: Stieltjes series}
\section{Stieltjes series}
\label{Sec:StieltjesSeries}

A function $F \colon \mathbb{C} \to \mathbb{C}$ is called a
\emph{Stieltjes} function if it can be expressed by the Stieltjes integral
\begin{equation}
  \label{StieFunF}
F (z) \; = \; \int_{0}^{\infty} \,
\frac{\mathrm{d} \Phi (t)} {z+t} \, ,
\qquad \vert \arg (z) \vert < \pi \, .
\end{equation}
Here, $\Phi (t)$ is a bounded, nondecreasing function taking infinitely
many different values on the interval $0 \le t < \infty$. Moreover, the
moment integrals
\begin{equation}
  \label{StieMom}
\mu_n \; = \; \int_{0}^{\infty} \, t^n \, \mathrm{d} \Phi (t) \, ,
\qquad n \in \mathbb{N}_0 \, ,
\end{equation}
must be positive and finite for all finite values of $n$.

Detailed discussions of Stieltjes series and their special role in the
theory of summability can be found in the books by Bender and Orszag
\cite[Chapter 8.6]{Bender/Orszag/1978} or Baker and Graves-Morris
\cite[Chapter 5]{Baker/Graves-Morris/1996}. In the case of divergent
Stieltjes series, it can be shown rigorously that the Pad\'{e}
approximants $[n+j/n]$ with fixed $j \ge -1$ converge to a uniquely
determined Stieltjes function as $n \to \infty$.

An inverse power series representation for such a Stieltjes function $F
(z)$ can be derived easily. We insert the geometric series
$\sum_{\nu=0}^{\infty} (-t)^{\nu}/z^{\nu+1} = 1/(z+t)$, which converges
for $\vert t/z \vert < 1$, into the integral representation
(\ref{StieFunF}) and -- ignoring all questions of legitimacy and
convergence -- integrate term-wise from $0$ to $\infty$ using
(\ref{StieMom}).  Thus, a Stieltjes function $F (z)$ can at least
formally be represented by its Stieltjes series
\begin{equation}
  \label{StieSerF}
  F (z) \; = \;
  \sum_{\nu=0}^{\infty} \, \frac{(-1)^{\nu} \, \mu_{\nu}}{z^{\nu+1}} \, .
\end{equation}

A factorial series for a a Stieltjes function $F (z)$ can also be derived
quite easily. For that purpose, we use the \emph{convergent} factorial
series (\ref{WaringFormula}) with $w = -t$ in the integral representation
(\ref{StieFunF}) and interchange integration and summation:
\begin{align}
  F (z) & \; = \; \int_{0}^{\infty} \, \sum_{n=0}^{\infty} \,
  \frac{(-t)_n}{(z)_{n+1}} \, \mathrm{d} \Phi (t)
  \\
  \label{StierSer_GenMomInt}
  & \; = \; \sum_{n=0}^{\infty} \, \frac{1}{(z)_{n+1}} \,
  \int_{0}^{\infty} \, (-t)_n \, \mathrm{d} \Phi (t) \, .
\end{align}
If we now expand $(-t)_n$ via (\ref{St1_GenFun1}) and do the resulting
moment integrals according to (\ref{StieMom}), we obtain:
\begin{align}
  \label{StieFun_MonIntSer}
  F (z) & \; = \; \sum_{n=0}^{\infty} \, \frac{(-1)^n}{(z)_{n+1}} \,
  \sum_{\nu=0}^{n} \, \textbf{S}^{(1)} (n, \nu) \, \int_{0}^{\infty} \,
  t^{\nu} \, \mathrm{d} \Phi (t)
  \\
  \label{StieFun_FactSer}
  & \; = \; \sum_{n=0}^{\infty} \, \frac{(-1)^n}{(z)_{n+1}} \,
  \sum_{\nu=0}^{n} \, \textbf{S}^{(1)} (n, \nu) \, \mu_{\nu} \, .
\end{align}
Obviously, the factorial series (\ref{StieFun_FactSer}) for a Stieltjes
function is a special case of the more general result
(\ref{InvPowSer->FactSer}).

Again, cancellation is the reason why we can hope that the factorial
series (\ref{StieFun_FactSer}) converges and is computationally useful
even if the corresponding Stieltjes series (\ref{StieSerF}) diverges. The
power $t^{n}$ in the moment integral (\ref{StieMom}) is positive, but the
Pochhammer symbol $(-t)_n$ in the generalized moment integrals
$\int_{0}^{\infty} (-t)_n \mathrm{d} \Phi (t)$ in
(\ref{StierSer_GenMomInt}) has zeros for $t = 0, 1, \dots n-1$ and
alternates in sign. This leads to a substantial cancellation.

\typeout{==> Section: Summation of the asymptotic series for the
  exponential integral}
\section{Summation of the asymptotic series for the exponential integral}
\label{Sec:SumAsySer_ExpInt}

The practical usefulness of the transformation of a factorially divergent
inverse power series to a convergent factorial series via
(\ref{InvPowSer->FactSer}) or via (\ref{StieFun_FactSer}) in the case of
a Stieltjes series can be demonstrated by means of the exponential
integral \cite[Eq.\ (5.1.1)]{Abramowitz/Stegun/1972}
\begin{equation}
  \label{Def_ExpInt}
  E_1 (z) \; = \;
  \int_{z}^{\infty} \, \frac{\exp (-t) \mathrm{d}t}{t} \, ,
\end{equation}
which possesses the following asymptotic expansion as $z \to \infty$
\cite[Eq.\ (5.1.51)]{Abramowitz/Stegun/1972},
\begin{equation}
  \label{ExpInt_AsySer}
  z \, \mathrm{e}^z \, E_1 (z) \; \sim \;
  \sum_{m=0}^{\infty} \, \frac{(-1)^m m!}{z^m}
  \; = \; {}_2 F_0 (1, 1; -1/z) \, , \qquad z \to \infty \, .
\end{equation}
The exponential integral can also be expressed as a Stieltjes integral
\cite[Eq.\ (5.1.28)]{Abramowitz/Stegun/1972}:
\begin{equation}
  \label{ExpInt_StieInt}
  \mathrm{e}^z \, E_1 (z) \; = \;
  \int_{0}^{\infty} \, \frac{\exp (-t) \mathrm{d}t}{z+t} \; = \;
  \frac{1}{z} \,
  \int_{0}^{\infty} \, \frac{\exp (-t) \mathrm{d}t}{1+t/z} \, .
\end{equation}
Thus, the divergent inverse power series (\ref{ExpInt_AsySer}) is a
Stieltjes series.

If we now combine (\ref{ExpInt_AsySer}) with (\ref{StieFun_MonIntSer})
and (\ref{StieFun_FactSer}), we obtain:
\begin{equation}
  \label{ExpInt_FactSer}
  \mathrm{e}^{z} \, E_{1} (z) \; = \; 
  \sum_{n=0}^{\infty} \, \frac{1}{(z)_{n+1}} \,
  \int_{0}^{\infty} \, (-t)_n \, \mathrm{e}^{-t} \mathrm{d} t 
  \; = \; \sum_{n=0}^{\infty} \, \frac{(-1)^n}{(z)_{n+1}} \,
  \sum_{\nu=0}^{n} \, \textbf{S}^{(1)} (n, \nu) \, \nu! \, .  
\end{equation}

As discussed in Sections \ref{Sec:InvPowSer<->FactSer} and
\ref{Sec:StieltjesSeries}, the transformation of an inverse power series
to a factorial series only produces a numerically useful result if a
substantial cancellation takes place in the inner sums in
(\ref{InvPowSer->FactSer}) or (\ref{StieFun_FactSer}). For that purpose,
the inverse power series coefficients $(-1)^{n} n!$ and the factorial
series coefficients $(-1)^{n} \sum_{\nu=0}^{n} \textbf{S}^{(1)} (n, \nu)
\nu!$ are displayed in Table \ref{Tab_6_1}.

\begin{table}[h]
\begin{center}
  \caption{Leading coefficients of the inverse power and the factorial
    series of the exponential integral.}
  \label{Tab_6_1} 
\begin{tabular*}{0.75\textwidth}{@{\extracolsep{\fill}}rrr}%
  \\
  \hline \hline %
  \multicolumn{1}{r}{$n$} &
  \multicolumn{1}{r}{$(-1)^{n} n!$} &
  \multicolumn{1}{r}{$(-1)^{n} 
  \sum_{\nu=0}^{n} \textbf{S}^{(1)} (n, \nu) \nu!$}%
  \rule[-7pt]{0pt}{19pt} \\
  \hline
  0   & 1 & 1 \\
  1   & -1 & -1 \\
  2   & 2 & 1 \\
  3   & -6 & -2 \\
  4   & 24 & 4 \\
  5   & -120 & -14 \\
  6   & 720 & 38 \\
  7   & -5~040 & -216 \\
  8   & 40~320 & 600 \\
  9   & -362~880 & -6~240 \\
  10  & 3~628~800 & 9~552 \\
  11  & -39~916~800 & -319~296 \\
  12  & 479~001~600 & -519~312 \\
  13  & -6~227~020~800 & -28~108~560 \\
  14  & 87~178~291~200 & -176~474~352 \\
  \hline \hline %
\end{tabular*}
\end{center}
\end{table}

The numbers displayed in Table \ref{Tab_6_1} show that there is indeed a
substantial amount of cancellation in the inner sums in
(\ref{ExpInt_FactSer}): The coefficients of the inverse power series grow
much more rapidly in magnitude than the coefficients of the factorial
series. This is enough to produce a convergent result. For example, the
first 15 coefficients of the factorial series produce for $z=5$ the
following result:
\begin{equation}
\frac{\displaystyle\sum_{n=0}^{14} \, \frac{(-1)^n}{(5)_{n+1}} \,
  \sum_{\nu=0}^{n} \, \textbf{S}^{(1)} (n, \nu) \, \nu!}
{\displaystyle\exp (5) E_{1} (5)} \; = \; 1.000~000~764
\end{equation}

\typeout{==> Section: Conversion of a power series to a factorial series}
\section{Conversion of a power series to a factorial series}
\label{Sec:PowSer2FactSer}

Assume that a function $f \colon \mathbb{C} \to \mathbb{C}$ possesses a
power series $f (z) = \sum_{n=0}^{\infty} \gamma_{n} z^{n}$. For a
transformation of this power series to a factorial series, we express it
as an inverse power series in $1/z$:
\begin{equation}
  \label{f_TransInvPow}
  f (z) \; = \; \frac{1}{z} \, \sum_{n=0}^{\infty} \, 
  \frac{\gamma_{n}}{(1/z)^{n+1}} \, .
\end{equation}
If we now use (\ref{InvPowSer->FactSer}), we obtain a factorial series in
$1/z$:
\begin{equation}
  \label{FactSer_f_TransInvPow_1}
  \sum_{n=0}^{\infty} \, \gamma_{n} \, z^{n} 
      \; = \; \frac{1}{z} \, 
              \sum_{m=0}^{\infty} \, \frac{(-1)^{m}}{(1/z)_{m+1}} \, 
              \sum_{\mu=0}^{m} \, (-1)^{\mu} \, 
              \textbf{S}^{(1)} (m, \mu) \, \gamma_{\mu} \, .
\end{equation}
An equivalent factorial series was considered by Ramis and Thomann
\cite[p.\ 20]{Ramis/Thomann/1981}. Thomann (\cite[p.\ 526]{Thomann/1990}
and \cite[Section 5.3]{Thomann/1995}) considered similar expansions in
terms of generalized factorial series, where the Pochhammer symbols are
replaced by products $z(z+t_{1}) \dots (z+t_{n})$. It is clear that
generalized factorial series are at least potentially more powerful than
their ordinary counterparts. However, it is not \emph{a priori} clear how
the parameters $\{ t_n \}_{n=1}^{\infty}$ should be chosen. In addition,
Thomann's formulas contain so-called generalized Stirling numbers instead
of the ordinary Stirling numbers of the first kind. Therefore, it is not
immediately obvious whether Thomann's generalized transformation formula
is really more useful than (\ref{FactSer_f_TransInvPow_1}).

Further manipulations of the Pochhammer symbol in
(\ref{FactSer_f_TransInvPow_1}) are possible:
\begin{equation}
  \label{VarSub_PoweSer->FactSer}
  \frac{1}{(1/z)_{m+1}} 
  \; = \; \frac{1}{\prod_{k=0}^{m} [k+1/z]}
  \; = \; \frac{z^{m+1}}{\prod_{k=1}^{m} k[z+1/k]}
   \; = \; \frac{z}{m!} \, \prod_{k=1}^{m} \, \frac{z}{z+1/k} \, .
\end{equation}
Inserting this into (\ref{FactSer_f_TransInvPow_1}) yields:
\begin{equation}
  \label{PoweSer->FactSer}
  \sum_{\nu=0}^{\infty} \, \gamma_{\nu} \, z^{\nu} 
  \; = \; \sum_{m=0}^{\infty} \, \frac{(-1)^{m}}{m!} \, 
          \prod_{k=1}^{m} \, \frac{z}{z+1/k} 
          \sum_{\mu=0}^{m} \, (-1)^{\mu} \, 
          \textbf{S}^{(1)} (m, \mu) \, \gamma_{\mu} \, .
\end{equation}
If the power series coefficients $\gamma_{\nu}$ have strictly alternating
signs, we can expect cancellation in the inner sum involving the Stirling
numbers. Moreover, for $z > 0$ we have $z/(z+1/k) < 1$ for $k \in
\mathbb{N}$. Since the transformation (\ref{VarSub_PoweSer->FactSer})
also produces the factorial $1/m!$, we can expect that the transformation
(\ref{PoweSer->FactSer}) produces a convergent and numerically useful
result even if the coefficients $\gamma_{\nu}$ diverge factorially in
magnitude.

\typeout{==> Section: The quartic anharmonic oscillator}
\section{The quartic anharmonic oscillator}
\label{Sec:TheQuarticAnharmonicOscillator}

In their seminal articles \cite{Bender/Wu/1969,Bender/Wu/1971}, Bender
and Wu showed that the Rayleigh-Schr\"{o}dinger perturbation expansions
for the energy eigenvalues $E^{(m)} (\beta)$ of the anharmonic
oscillators defined by the Hamiltonians
\begin{equation}
  \label{Ham_AHM}
  \hat{H}^{(m)} (\beta) \; = \; \hat{p}^{2} + \hat{x}^{2} +
  \beta \hat{x}^{2m} \, , \qquad m = 2, 3, 4, \ldots \, , \qquad 
  \hat{p} \; = \; - \mathrm{i} \frac{\mathrm{d}}{\mathrm{d} x} \, ,
\end{equation}
diverge quite violently for every nonzero coupling constant $\beta$
(here, the same notation as in
\cite{Weniger/Cizek/Vinette/1993,Weniger/1996c,Weniger/1996e} is used).
Later, perturbation expansions with a similar type of divergence were
discovered in the case of other quantum mechanical systems, which
ultimately created a new sub-discipline of theoretical physics called
\emph{large order perturbation theory} (see for example the book by Le
Guillou and Zinn-Justin \cite{LeGuillou/Zinn-Justin/1990} and the
articles reprinted there).

In the following years, a lot of work has been done on the summation of
divergent perturbation expansions (a mathematically oriented overview can
be found in \cite{Simon/1991}). In particular the quartic anharmonic
oscillator with $m=2$ in (\ref{Ham_AHM}) has become a very popular
computational laboratory for theoretical physicists. In spite of its
simplicity, the quartic anharmonic oscillator leads to challenging
computational and conceptual problems, as documented in countless
articles (far too many to be cited here).

In this article, I am exclusively interested in the summation of the
divergent Rayleigh-Schr\"{o}dinger perturbation expansion
\begin{equation}
  \label{Eq_beta}
  E^{(2)}(\beta) \; = \;
  \sum_{n=0}^{\infty} \, b_{n}^{(2)} \, \beta^{n}
\end{equation}
for the ground state energy of the quartic anharmonic oscillator by
transforming it to a factorial series via (\ref{PoweSer->FactSer}).
Long, but nevertheless incomplete lists of references dealing with other
approaches for the summation of the divergent perturbation expansions of
the anharmonic oscillators can be found in
\cite{Weniger/Cizek/Vinette/1993,Weniger/1996c,Weniger/1996e} or also in
\cite[Kap.\ 10]{Weniger/1994b}.

If the convention (\ref{Ham_AHM}) for the Hamiltonian is used, the
coefficients $b_{n}^{(2)}$ possess in the case of large indices $n$ the
following leading order asymptotics (see for example \cite[Eq.\
(2.3)]{Weniger/1996c}):
\begin{equation}
  \label{bAs2}
  b_{n}^{(2)} \; \sim \; (-1)^{n+1} \,
  \frac{(24)^{1/2}}{\pi^{3/2}} \, \Gamma(n+1/2) \, (3/2)^{n} \, ,
  \qquad n \to \infty \, .
\end{equation}
This asymptotic estimate shows that the perturbation series
(\ref{Eq_beta}) diverges for all $\beta \ne 0$ like the generalized
hypergeometric series ${}_2 F_0 (1/2, 1; - 3 \beta / 2) =
\sum_{m=0}^{\infty} (1/2)_m (- 3 \beta / 2)^m$ \cite[Eq.\
(1.10)]{Weniger/1990}.

It was shown rigorously bu Simon \cite[Theorem IV.2.1]{Simon/1970a} that
the perturbation expansion
\begin{equation}
  \label{Delta_Eq_beta}
  \Delta E^{(2)}(\beta) \; = \;
  \sum_{n=0}^{\infty} \, b_{n+1}^{(2)} \, \beta^{n}
\end{equation}
for the energy shift defined by $E^{(2)}(\beta) = b_{0}^{(2)} + \beta
\Delta E^{(2)}(\beta) = 1 + \beta \Delta E^{(2)}(\beta)$ is a Stieltjes
series. This is a highly advantageous feature. As discussed in Section
\ref{Sec:StieltjesSeries}, this implies that the perturbation series
(\ref{Delta_Eq_beta}) corresponds to a uniquely defined Stieltjes
function since it is Pad\'{e} summable. Moreover, the asymptotic estimate
(\ref{bAs2}) implies that the terms of the perturbation expansion
(\ref{Delta_Eq_beta}) have for $\beta > 0$ strictly alternating signs
(see also \cite[Tabelle 10-1]{Weniger/1994b}), which is advantageous if
we want to sum it with the help of (\ref{PoweSer->FactSer}).

If we transform the perturbation series (\ref{Delta_Eq_beta}) for the
energy shift with the help of (\ref{PoweSer->FactSer}), we obtain the
following expansion for the ground state energy of the quartic
oscillator:
\begin{equation}
  \label{Eq_beta->FactSer}
  E^{(2)}(\beta) \; = \; 1 \, + \, \beta \, 
    \sum_{m=0}^{\infty} \, \frac{(-1)^{m}}{m!} \, 
      \prod_{k=1}^{m} \, \frac{\beta}{\beta+1/k} 
        \sum_{\mu=0}^{m} \, (-1)^{\mu} \, 
          \textbf{S}^{(1)} (m, \mu) \, b_{\mu+1}^{({2})} \, .
\end{equation}
The first 34 terms of the infinite series on the right-hand side yield
for $\beta = 1/5$ the energy $E_{\mathrm{FS}}^{(2)} (1/5) = 1.118~305
\dots$, which is less accurate than the energy $E_{\mathrm{PA}}^{(2)}
(1/5) \; = \; 1.118~292~654~373 \dots$ obtained by computing the Pad\'{e}
approximants $[17/17]$ from the first 34 terms of the perturbation
expansion (\ref{Delta_Eq_beta}) for the energy shift $\Delta
E^{(2)}(1/5)$. These approximations can be compared to the ``exact''
energy $E_{\mathrm{exact}}^{(2)} (1/5) \; = \; 1.118~292~654~367~039~154
\dots$ obtained by a very sophisticated summation calculation
\cite[Tabelle 10-9]{Weniger/1994b}.

Thus, the truncated series expansion (\ref{Eq_beta->FactSer}), which does
a \emph{linear} transformation of the perturbation series coefficients
$b_{1}$, $b_{2}$, $\dots$, $b_{34}$, is less efficient than the highly
\emph{nonlinear} Pad\'{e} approximant using the same number of
coefficients $b_{\mu}^{(2)}$. Nevertheless, improvements are possible if
we use the integral representation (\ref{FS_IntRep}) for the evaluation
of the corresponding factorial series $\Omega (z) = \sum_{n=0}^{\infty}
a_{n} n!/(z)_{n+1}$.

A direct use of a truncation of the power series $\varphi (t) =
\sum_{n=0}^{\infty} a_{n} (1-t)^{n}$ defined by (\ref{FS_IntRep_b}) in
the integral representation (\ref{FS_IntRep_a}) does not lead to an
improvement since integration is linear. It is, however, possible to
replace the truncated power series for $\varphi (t)$ by a Pad\'{e}
approximant in $1-t$ and to evaluate the resulting expression by
numerical quadrature. The use of the Pad\'{e} approximant $[17/17]$ to
$\varphi (t)$ in (\ref{FS_IntRep}) yields $E_{\mathrm{IntFS}} (1/5) =
1.118~292~654~369 \dots$, which is better than the direct Pad\'{e}
summation result $E_{\mathrm{PA}}^{(2)} (1/5) \; = \; 1.118~292~654~373
\dots$, but less accurate than $E_{\mathrm{BP}} (1/5) \; =
\;1.118~292~654~367~039~152 \dots$ obtained by doing a so-called
Borel-Pad\'{e} transformation originally introduced in
\cite{Graffi/Grecchi/Simon/1970}.

\typeout{==> Section: Outlook}
\section{Outlook}
\label{Sec:Outlook}

A single article cannot provide an exhaustive treatment of the numerical
utilization of factorial series for the summation of divergent (inverse)
power series. Because of space limitations, many interesting or
potentially useful aspect of the theory of factorial series were treated
only superficially or even completely ignored. Nevertheless, I hope that
this article will inspire others.

Inner sums of the type of $(-1)^{m} \sum_{\mu=0}^{m} (-1)^{\mu}
\textbf{S}^{(1)} (m, \mu) c_{\mu}$ occurring on the right-hand side of
(\ref{InvPowSer->FactSer}) are the key quantities of this article.
Recurrence formulas or even alternative closed form expressions would
obviously be desirable. In the case of Stieltjes series, a recursive
scheme for the computation of the generalized moment integrals
$\int_{0}^{\infty} \, (-t)_n \, \mathrm{d} \Phi (t)$ in
(\ref{StierSer_GenMomInt}) can be derived.

The numerical examples of this article also raise questions. It is an
obvious question whether and how well the convergence of factorial series
can be accelerated by sequence transformations. When I looked at the
convergence of the factorial series (\ref{ExpInt_FactSer}) for $E_{1}
(z)$, I applied sequence transformations to speed up its convergence.
However, convergence was not improved substantially by the
transformations I used.  At the moment, it is unclear whether this is a
specific feature of the factorial series (\ref{ExpInt_FactSer}), or
whether we face a problem of a more general nature. It could be that the
convergence of factorial series can only be accelerated effectively if
other, specially designed sequence transformations are used. This should
be investigated.

It should be worthwhile to investigate whether the transformation formula
(\ref{PoweSer->FactSer}) can also be used profitably in the case of
\emph{convergent} power series as a convergence acceleration tool. One
can also hope that in the case of sufficiently simple power series
coefficients $\gamma_{n}$ explicit expressions for the inner sum
$(-1)^{m} \sum_{\mu=0}^{m} (-1)^{\mu} \textbf{S}^{(1)} (m, \mu)
\gamma_{\mu}$ can be found, which would yield new explicit expressions in
terms of factorial series for functions defined by power series.

In Section \ref{Sec:TheQuarticAnharmonicOscillator}, the truncated power
series $\varphi (t) = \sum_{n=0}^{\infty} a_{n} (1-t)^{n}$ defined by
(\ref{FS_IntRep_b}) was converted to a Pad\'{e} approximant in $1-t$,
which was inserted into the integral representation (\ref{FS_IntRep}).
One should investigate whether other sequence transformations produce
better approximations to $\varphi (t)$ than Pad\'{e} approximants.

\appendix

\typeout{==> Appendix A: Stirling numbers}
\section{Stirling numbers}
\label{App:StirlingNumbers}
 
With respect to notation, the theory of Stirling numbers is a mess. This
is partly due to the fact that numerous different symbols are used in the
literature (a discussion of the various notations can be found in
\cite[p.\ 822]{Abramowitz/Stegun/1972}). To make things worse, different
and incompatible notations are used for factorial expressions. In special
function theory, Pochhammer symbols $(z)_{n} = z(z+1) \dots (z+n-1) =
\Gamma (z+n)/\Gamma (z)$ are consistently used, but in combinatorics, it
is more common to use instead falling factorials $z(z-1)\dots (z-n+1) =
\Gamma(z+1)/ \Gamma(z-n+1)$. Unfortunately, in the literature of
combinatorics falling factorials are often denoted by the symbol
$(z)_{n}$ normally reserved in special function theory for Pochhammer
symbols.

The Stirling numbers $\textbf{S}^{(1)} (n, \nu)$ of the first kind are
the polynomial coefficients of a Pochhammer symbol $(z-n+1)_n = z(z-1)
\dots (z-n+1) = \Gamma (z+1)/\Gamma (z-n+1)$ (see for example \cite[Eq.\
(1) on p.\ 56]{Srivastava/Choi/2001}):
\begin{equation}
  \label{St1_GenFun1}
(z-n+1)_n \; = \; (-1)^n (-z)_n \; = \;
\sum_{\nu=0}^{n} \, \textbf{S}^{(1)} (n, \nu) \, z^{\nu} \, ,
\qquad n \in \mathbb{N}_0 \, .
\end{equation}
If we use $(z-n+1)_n = (-1)^{n} (-z)_{n}$ and replace $z$ by $-z$, we
obtain:
\begin{equation}
  \label{St1_GenFun2}
(z)_n \; = \; (-1)^n \, \sum_{\nu=0}^{n} \, (-1)^{\nu} \,
\textbf{S}^{(1)} (n, \nu) \, z^{\nu} \, ,
\qquad n \in \mathbb{N}_0 \, .
\end{equation}
If $z > 0$ holds, the coefficients of all powers $z^{\nu}$ with $0 \le
\nu \le n$ in the expansion of $(z)_{n}$ are either zero or positive.
Thus, (\ref{St1_GenFun2}) implies
\begin{equation}
  \label{St1_SignConv}
  (-1)^{n-\nu} \, \textbf{S}^{(1)} (n, \nu)  \; = \;
  \bigl\vert \textbf{S}^{(1)} (n, \nu) \bigr\vert \, .
\end{equation}

The Stirling numbers $\textbf{S}^{(2)} (n, \nu)$ of the second kind are
usually defined as follows (see for example \cite[Eq.\ (14) on p.\
58]{Srivastava/Choi/2001}):
\begin{equation}
  \label{St2_GenFun1}
  z^n \; = \;
  \sum_{\nu=0}^{n} \, \textbf{S}^{(2)} (n, \nu) \, (z-\nu+1)_\nu \, ,
  \qquad n \in \mathbb{N}_0 \, . 
\end{equation}

The Stirling numbers of the first and second kind correspond to
triangular matrices that transform the polynomial sequences $\{
(z-n+1)_{n} \}_{n=0}^{\infty}$ and $\{ z^{n} \}_{n=0}^{\infty}$ into each
other. Since these transformations must be invertible, orthogonality
relationships exist.

If we replace in (\ref{St1_GenFun1}) the powers on the right-hand side by
Pochhammer symbols according to (\ref{St2_GenFun1}), we obtain:
\begin{equation}
(z-n+1)_n \; = \; \sum_{\nu=0}^{n} \, \textbf{S}^{(1)} (n, \nu) \,
\sum_{k=0}^{\nu} \, \textbf{S}^{(2)} (\nu, k) \, (z-k+1)_{k} \, ,
\qquad n \in \mathbb{N}_0 \, .
\end{equation}
By interchanging the order of the summations summations, this expression
can be rewritten as follows:
\begin{equation}
(z-n+1)_n \; = \; \sum_{k=0}^{n} \, (z-k+1)_k \,
\sum_{\nu=k}^{n} \, \textbf{S}^{(1)} (n, \nu) \,
\textbf{S}^{(2)} (\nu, k) \, , \qquad n \in \mathbb{N}_0 \, .
\end{equation}
Thus, we obtain the following well-known orthogonality relationship:
\begin{equation}
  \label{St1St2_OrthRel1}
\sum_{\nu=k}^{n} \, \textbf{S}^{(1)} (n, \nu) \,
\textbf{S}^{(2)} (\nu, k) \; = \; \delta_{n k} \, ,
\qquad k, n \in \mathbb{N}_0 \, .
\end{equation}

By replacing in the Pochhammer symbols on the right-hand side of
(\ref{St2_GenFun1}) by powers according to (\ref{St1_GenFun1}), we obtain
the following alternative orthogonality relationship:
\begin{equation}
  \label{St2St1_OrthRel1}
\sum_{\nu=k}^{n} \, \textbf{S}^{(2)} (n, \nu) \,
\textbf{S}^{(1)} (\nu, k) \; = \; \delta_{n k} \, ,
\qquad k, n \in \mathbb{N}_0 \, .
\end{equation} 

The Stirling numbers $\textbf{S}^{(2)} (n, \nu)$ of the second kind
possess the following infinite generating function (see for example
\cite[Eq.\ (16) on p.\ 58]{Srivastava/Choi/2001}):
\begin{equation}
  \label{GenFunSti2}
  \frac{1}{(1-t)(1-2t)\dots(1-kt)} \; = \; \sum_{\kappa=0}^{\infty} \,
  \textbf{S}^{(2)} (k+\kappa, k) \, t^{\kappa} \, , 
  \qquad k \in \mathbb{N} \, , \quad \vert t \vert < 1/k \, .
\end{equation}
The substitution $t = 1/z$ yields:
\begin{equation}
\frac{1}{(z-k)_{k+1}} \; = \;
\sum_{\kappa=0}^{\infty} \, \textbf{S}^{(2)} (k+\kappa, k) \,
z^{-k-\kappa-1} \, , \qquad \vert z \vert > k \, .
\end{equation}
If we now use $(z-k)_{k+1} = (-1)^{k+1} (-z)_{k+1}$ and replace $z$ by
$-z$, we obtain the following inverse power series expansion of an
inverse Pochhammer symbol \cite[Eq.\ (9) on p.\ 68]{Nielsen/1965}:
\begin{equation}
  \label{Poch_expand_Pow}
\frac{1}{(z)_{k+1}} \; = \; \sum_{\kappa=0}^{\infty} \,
\frac{(-1)^{\kappa} \, \textbf{S}^{(2)} (k+\kappa, k)}{z^{k+\kappa+1}}
\, , \qquad k \in \mathbb{N}_{0} \, , \quad \vert z \vert > k \, .
\end{equation}

A convenient starting point for the derivation of a factorial series for
an inverse power $1/z^{k+1}$ is the following factorial series \cite[Eq.\
(3) on p.\ 77]{Nielsen/1965}:
\begin{equation}
  \label{WaringFormula}
\frac{1}{z-w} \; = \; \sum_{n=0}^{\infty} \,
\frac{(w)_n}{(z)_{n+1}} \, , \qquad \Re (z-w) > 0 \, .
\end{equation}
Next, we apply $\mathrm{d}^k [1/(z-w)]/\mathrm{d} w^k = k!/(z-w)^{k+1}$
with $k \in \mathbb{N}_0$ to (\ref{WaringFormula}), which yields:
\begin{equation}
  \label{DiffWaringFormula}
  \frac{k!}{(z-w)^{k+1}} \; = \; \frac{\mathrm{d}^k}{\mathrm{d} w^k}
  \, \sum_{n=0}^{\infty} \, \frac{(w)_n}{(z)_{n+1}} \, .
\end{equation}
The Pochhammer symbol $(w)_n$ can be expanded with the help of
(\ref{St1_GenFun1}). Then, the differentiations can be done in closed
form and we obtain -- after setting $w=0$ -- the following factorial
series for an inverse power \cite[Eq.\ (6) on p.\ 78]{Nielsen/1965}:
\begin{equation}
  \label{InvPow_expand_InvPoch}
\frac{1}{z^{k+1}} \; = \; \sum_{\kappa=0}^{\infty} \,
\frac{(-1)^{\kappa} \, \textbf{S}^{(1)} (k+\kappa, k)}{(z)_{k+\kappa+1}}
\, , \qquad k \in \mathbb{N}_0 \, .
\end{equation}

\typeout{==> Appendix B: Sequence inversion by triangular orthogonal
  matrices}
\section{Sequence inversion by triangular orthogonal matrices}
\label{App:SequenceInversionByTriangularMatrices}

Let us assume that we have two matrices $\mathbf{A} = \{ A_{m
  n} \}_{m, n \ge 0}$ and $\mathbf{B} = \{ B_{m n} \}_{m, n \ge 0}$,
which are triangular, 
\begin{equation}
  \label{A_B_trinang}
  A_{m n} \; = \; B_{m n} \; = \; 0 \, , \qquad n > m \, ,
\end{equation}
and which satisfy the orthogonality relationships
\begin{align}
  \label{Orth_B_A}
  \sum_{r=k}^{n} \, B_{n r} \, A_{r k} & \; = \; \delta_{n k} \, ,
  \\
  \label{Orth_A_B}
  \sum_{r=k}^{n} \, A_{n r} \, B_{r k} & \; = \; \delta_{n k} \, .  
\end{align}
The orthogonality relationships (\ref{St1St2_OrthRel1}) and
(\ref{St2St1_OrthRel1}) involving the Stirling numbers of the first and
second kind are an example of such a pair of orthogonality relationships.

We also assume that there are two sequences $\{ x_n \}_{n=0}^{\infty}$
and $\{ y_n \}_{n=0}^{\infty}$ whose elements are connected by finite
linear combinations
\begin{equation}
  \label{x_n->y_n}
  y_{r} \; = \; \sum_{k=0}^{r} \, A_{r k} \, x_{k} \, , 
  \qquad r \in \mathbb{N}_{0} \, ,
\end{equation}
and two sequences $\{ u_n \}_{n=0}^{\infty}$ and $\{
w_n \}_{n=0}^{\infty}$, whose elements are connected by infinite
series expansions
\begin{equation}
  \label{w_n->u_n}
  u_{r} \; = \; \sum_{n=r}^{\infty} \, A_{n r} \, w_{n} \, ,
  \qquad r \in \mathbb{N}_{0} \, . 
\end{equation}

In the book of Charalambides \cite[Example 8.2 on pp.\ 284 -
285]{Charalambides/2002} it is shown that the finite linear combination
(\ref{x_n->y_n}) then possesses the inverse relation
\begin{equation}
  \label{y_n->x_n}
  x_{n} \; = \; \sum_{k=0}^{n} \, B_{n k} \, y_{k} \, , 
  \qquad n \in \mathbb{N}_{0} \, , 
\end{equation}
and that the infinite series expansion (\ref{w_n->u_n}) then possesses
the inverse relation
\begin{equation}
  \label{u_n->w_n}
  w_{k} \; = \; \sum_{r=k}^{\infty} \, B_{r k} \, u_{r} \, , 
  \qquad k \in \mathbb{N}_{0} \, .  
\end{equation} 

Obviously, the two finite linear combinations (\ref{x_n->y_n}) and
(\ref{y_n->x_n}) generalize the finite generating functions
(\ref{St1_GenFun1}) and (\ref{St2_GenFun1}) of the Stirling numbers, and
the two infinite series expansions (\ref{w_n->u_n}) and (\ref{u_n->w_n})
generalize the series expansions (\ref{Poch_expand_Pow}) and
(\ref{InvPow_expand_InvPoch}) connecting inverse powers and inverse
Pochhammer symbols.

Further generalizations are possible. Let us now consider the following,
in general \emph{formal} infinite series expansion:
\begin{equation}
  \label{Def_F}
  F \; = \; \sum_{n=0}^{\infty} \, \eta_{n} \, y_{n} \, .
\end{equation}
With the help of (\ref{x_n->y_n}), the elements of the sequence $\{ y_n
\}_{n=0}^{\infty}$ can be replaced by the elements of the sequence $\{
x_n \}_{n=0}^{\infty}$:
\begin{equation}
  \label{F_y->x_1}
  F \; = \; \sum_{n=0}^{\infty} \, \eta_{n} \, 
  \sum_{k=0}^{n} \, A_{n k} \, x_{k} \, .
\end{equation}
The order of the summations of this expansion can be rearranged:
\begin{align}
  \label{F_y->x_2}
  F & \; = \; x_{0} \, \bigl[ \eta_{0} \, A_{0 0} + \eta_{1} \, A_{1 0} +
  \eta_{2} \, A_{2 0} + \dots \bigr] \, + \, x_{1} \, \bigl[ \eta_{1} \,
  A_{1 1} + \eta_{2} \, A_{2 1} + \eta_{3} \, A_{3 1} + \dots \bigr]
  \notag \\
  & \phantom{\; = \;} \, + \, x_{2} \, \bigl[ \eta_{2} \, A_{2 2} +
  \eta_{3} \, A_{3 2} + \eta_{4} \, A_{4 2} + \dots \bigr] \, + \, x_{3}
  \, \bigl[ \eta_{3} \, A_{3 3} + \eta_{4} \, A_{4 3} + \eta_{5} \, A_{5
    3} + \dots \bigr] \, + \, \dots
  \\
  \label{F_y->x_3}
  & \; = \; \sum_{n=0}^{\infty} \, x_{n} \, \sum_{k=n}^{\infty} \, A_{k
    n} \, \eta_{k} \, .
\end{align}
Thus, the substitution $\{ y_n \}_{n=0}^{\infty} \to \{ x_n
\}_{n=0}^{\infty}$ in (\ref{Def_F}) produces an alternative expansion for
$F$ having the same general structure:
\begin{equation}
  \label{F_y->x_4}
  F \; = \; \sum_{n=0}^{\infty} \, \xi_{n} \, x_{n} \, ,
  \qquad
  \xi_{n} \; = \; \sum_{k=n}^{\infty} \, A_{k n} \, \eta_{k}
  \; = \; \sum_{\nu=0}^{\infty} \, A_{n+\nu, n} \, \eta_{n+\nu} \, .  
\end{equation}

Unknowingly, I had used transformation formulas for expansion
coefficients of the kind of (\ref{F_y->x_4}) already in
\cite{Weniger/2008} when I studied the transformation of Laguerre series
to power series.

Next, we consider the following, in general \emph{formal} infinite series
expansion:
\begin{equation}
  \label{Def_G}
  G \; = \; \sum_{n=0}^{\infty} \, \upsilon_{n} \, u_{n} \, .
\end{equation}
With the help of (\ref{w_n->u_n}), the elements of $\{ u_n
\}_{n=0}^{\infty}$ can be replaced by the elements of $\{ w_n
\}_{n=0}^{\infty}$:
\begin{equation}
  \label{G_u->w_1}
  G \; = \; \sum_{n=0}^{\infty} \, \upsilon_{n} \, 
  \sum_{k=n}^{\infty} \, A_{k n} \, w_{k} \, .
\end{equation}
The order of the summations of this expansion can be rearranged:
\begin{align}
  \label{G_u->w_2}
  G & \; = \; w_{0} \, \upsilon_{n} \, A_{0 0} \, + \, w_{1} \, \bigl[
  A_{1 0} \, \upsilon_{0} + A_{1 1} \, \upsilon_{1} \bigr] \, + \, w_{2}
  \, \bigl[ A_{2 0} \, \upsilon_{0} + A_{2 1} \, \upsilon_{1} + A_{2 2}
  \, \upsilon_{2} \bigr] \, + \, \dots
  \\
  \label{G_u->w_3}
  & \; = \; \sum_{n=0}^{\infty} \, w_{n} \, \sum_{k=0}^{n} \, A_{n k} \,
  \upsilon_{k} \, .
\end{align}
Thus, the substitution $\{ u_n \}_{n=0}^{\infty} \to \{ w_n
\}_{n=0}^{\infty}$ in (\ref{Def_G}) produces an alternative expansion of
the same general structure:
\begin{equation}
  \label{G_u->w_4}
  G \; = \; \sum_{n=0}^{\infty} \, \omega_{n} \, w_{n} \, ,
  \qquad 
  \omega_{n} \; = \; \sum_{k=0}^{n} \, A_{n k} \, \upsilon_{k} \, .
\end{equation}
These relationships for $G$ obviously generalize the transformation
formulas (\ref{InvPowSer->FactSer}) and (\ref{FactSer->InvPowSer}), which
transform factorial series and inverse power series into each other.

Let me emphasize once more that all operations considered in this Section
are purely \emph{formal}. Accordingly, we cannot tacitly assume that the
transformation formulas necessarily lead to convergent expansions. This
has to be checked explicitly in each case.

This Appendix was inspired by the orthogonality relationships
(\ref{St1St2_OrthRel1}) and (\ref{St2St1_OrthRel1}) involving the
Stirling numbers of the first and second kind. But many other
mathematical objects possess similar features. Recently, I had studied
the transformation of Laguerre expansions to power series expansions
\cite{Weniger/2008}. As is well known, generalized Laguerre polynomials
and powers are connected by finite sums of the type of (\ref{x_n->y_n})
and (\ref{y_n->x_n}) which can be inverted. Accordingly, these
transformation formulas possess certain orthogonality properties when
written in matrix form. Moreover, the coefficients of Laguerre series and
of power series are connected by infinite series expansions of the type
of (\ref{F_y->x_4}). I only understood the wider significance of these
features after \cite{Weniger/2008} was published when I studied the book
by Charalambides \cite{Charalambides/2002} more carefully.



\begin{thebibliography}{10}
\expandafter\ifx\csname url\endcsname\relax
  \def\url#1{\texttt{#1}}\fi
\expandafter\ifx\csname urlprefix\endcsname\relax\def\urlprefix{URL }\fi

\bibitem{Abramowitz/Stegun/1972}
M.~Abramowitz, I.~A. Stegun (eds.), Handbook of Mathematical Functions,
  National Bureau of Standards, Washington, D. C., 1972.

\bibitem{Baker/Graves-Morris/1996} G.~A. Baker, Jr., P.~Graves-Morris,
  {P}ad\'{e} Approximants, 2nd ed., Cambridge U. P., Cambridge, 1996.

\bibitem{Bender/Orszag/1978} C.~M. Bender, S.~A. Orszag, Advanced
  Mathematical Methods for Scientists and Engineers, McGraw-Hill, New
  York, 1978.

\bibitem{Bender/Wu/1969} C.~M. Bender, T.~T. Wu, Anharmonic oscillator,
  Phys. Rev. 184 (1969) 1231 -- 1260.

\bibitem{Bender/Wu/1971} C.~M. Bender, T.~T. Wu, Large-order behavior of
  perturbation theory, Phys. Rev.  Lett. 27 (1971) 461 -- 465.

\bibitem{Borel/1988} E.~Borel, Le\c{c}ons sur les S\'{e}ries Divergentes,
  2nd ed., \'{E}ditions Jacques Gabay, Paris, 1988. Originally published
  by Gautier-Villars, Paris, 1928. English translation by C. L.
  Critchfield and A. Vakar, Lectures on Divergent Series, Translation
  LA-6140-TR, Los Alamos Scientific Laboratory, Los Alamos, 1975.

\bibitem{Borghi/2007} R.~Borghi, Evaluation of diffraction catastrophes
  by using {W}eniger transformation, Opt. Lett. 32 (2007) 226 -- 228.

\bibitem{Borghi/2008c} R.~Borghi, Joint use of the {W}eniger
  transformation and hyperasymptotics for accurate asymptotic evaluations
  of a class of saddle-point integrals, Phys.  Rev. E 78 (2008) 026703--1
  -- 026703--11.

\bibitem{Borghi/2008b} R.~Borghi, On the numerical evaluation of cuspoid
  diffraction catastrophes, J.  Opt. Soc. Amer. A 25 (2008) 1682 -- 1690.

\bibitem{Borghi/2008a} R.~Borghi, Summing {P}auli asymptotic series to
  solve the wedge problem, J.  Opt. Soc. Amer. A 25 (2008) 211 -- 218.

\bibitem{Borghi/2009} R.~Borghi, Joint use of the {W}eniger
  transformation and hyperasymptotics for accurate asymptotic evaluations
  of a class of saddle-point integrals. {II}.  {H}igher-order
  transformations, Phys. Rev. E 80 (2009) 016704--1 -- 016704--15.

\bibitem{Borghi/Alonso/2007} R.~Borghi, M.~Alonso, Efficient evaluation
  of far-field asymptotic series, in: ICTON '07: 9th International
  Conference on Transparent Optical Networks, vol.~3, IEEE Xplore,
  Piscataway, NJ, 2007.

\bibitem{Borghi/Santarsiero/2003} R.~Borghi, M.~Santarsiero, Summing
  {L}ax series for nonparaxial beam propagation, Opt. Lett. 28 (2003) 774
  -- 776.

\bibitem{Brezinski/Matos/1996} C.~Brezinski, A.~C. Matos, A derivation of
  extrapolation algorithms based on error estimates, J. Comput. Appl.
  Math. 66 (1996) 5 -- 26.

\bibitem{Brezinski/RedivoZaglia/1994a} C.~Brezinski, M.~Redivo~Zaglia,
  {A} general extrapolation procedure revisited, Adv. Comput. Math. 2
  (1994) 461 -- 477.

\bibitem{Brezinski/RedivoZaglia/1994b} C.~Brezinski, M.~Redivo~Zaglia, On
  the kernel of sequence transformations, Appl. Numer. Math. 16 (1994)
  239 -- 244.

\bibitem{Caliceti/Meyer-Hermann/Ribeca/Surzhykov/Jentschura/2007}
  E.~Caliceti, M.~Meyer-Hermann, P.~Ribeca, A.~Surzhykov, U.~D.
  Jentschura, From useful algorithms for slowly convergent series to
  physical predictions based on divergent perturbative expansions, Phys.
  Rep. 446 (2007) 1 -- 96.

\bibitem{Charalambides/2002} C.~A. Charalambides, Enumerative
  Combinatorics, Chapman \& Hall, Boca Raton, 2002.

\bibitem{Cizek/Zamastil/Skala/2003} J.~{\v{C}\'{\i}\v{z}ek}, J.~Zamastil,
  L.~{Sk\'{a}la}, New summation technique for rapidly divergent
  perturbation series. {H}ydrogen atom in magnetic field, J. Math. Phys.
  44 (2003) 962 -- 968.

\bibitem{Delabaere/Rasoamanana/2007} E.~Delabaere, J.-M. Rasoamanana,
  Sommation effective d'une somme {B}orel par s\'{e}ries de factorielles,
  Annal. l'Inst. Fourier 57 (2007) 421 -- 456.

\bibitem{Graffi/Grecchi/Simon/1970} S.~Graffi, V.~Grecchi, B.~Simon,
  {B}orel summability: Application to the anharmonic oscillator, Phys.
  Lett. B 32 (1970) 631 -- 634.

\bibitem{Knopp/1964} K.~Knopp, Theorie und {A}nwendung der unendlichen
  {R}eihen, Springer-Verlag, Berlin, 1964.

\bibitem{Landau/1906} E.~Landau, {\"U}ber die {G}rundlagen der {T}heorie
  der {F}akult\"{a}tenreihen, Sitzungsb. K\"{o}nigl. Bay. Akad.
  Wissensch. M\"{u}nchen, math.-phys. Kl. 36 (1906) 151 -- 218.

\bibitem{LeGuillou/Zinn-Justin/1990} J.~C. Le~Guillou, J.~Zinn-Justin
  (eds.), Large-Order Behaviour of Perturbation Theory, North-Holland,
  Amsterdam, 1990.

\bibitem{Levin/1973} D.~Levin, Development of non-linear transformations
  for improving convergence of sequences, Int. J. Comput. Math. B 3
  (1973) 371 -- 388.

\bibitem{Li/Zang/Tian/2009} J.~Li, W.~Zang, J.~Tian, Simulation of
  {G}aussian laser beams and electron dynamics by {W}eniger
  transformation method, Opt. Expr. 17 (2009) 4959 -- 4969.

\bibitem{Li/Zang/Li/Tian/2009} J.-X. Li, W.~Zang, Y.-D. Li, J.~Tian,
  Acceleration of electrons by a tightly focused intense laser beam, Opt.
  Expr. 17 (2009) 11850 -- 11859.

\bibitem{Malgrange/1995} B.~Malgrange, Sommation des séries divergentes,
  Expo. Math. 13 (1995) 163 -- 222.

\bibitem{Milne-Thomson/1981} L.~M. Milne-Thomson, The Calculus of Finite
  Differences, Chelsea, New York, 1981. Originally published by
  Macmillan, London, 1933.

\bibitem{Nevanlinna/1919} F.~Nevanlinna, Zur {T}heorie asymptotischer
  {P}otenzreihen, Ann. Acad. Sci.  Fenn. Ser. A 12 (1919) 1 -- 81.

\bibitem{Nielsen/1965} N.~Nielsen, Die Gammafunktion, Chelsea, New York,
  1965. Originally published by Teubner, Leipzig and Berlin, 1906.

\bibitem{Noerlund/1926} N.-E. N\"{o}rlund, Le{\c c}ons sur les S\'eries
  d'Interpolation, Gautier-Villars, Paris, 1926.

\bibitem{Noerlund/1929} N.-E. N\"{o}rlund, Le{\c c}ons sur les
  \'{E}quations Lin\'{e}aires aux Diff\'{e}rences Finies,
  Gautier-Villars, Paris, 1929.

\bibitem{Noerlund/1954} N.~E. N\"{o}rlund, Vorlesungen \"{u}ber
  Differenzenrechnung, Chelsea, New York, 1954. Originally published by
  Springer-Verlag, Berlin, 1924.

\bibitem{OldeDaalhuis/2004b} A.~B. Olde~Daalhuis, Inverse
  factorial-series solutions of difference equations, Proc. Edinb. Math.
  Soc. 47 (2004) 421 -- 448.

\bibitem{Ramis/Thomann/1981} J.~P. Ramis, J.~Thomann, Some comments about
  the numerical utilization of factorial series, in: J.~Della~Dora,
  J.~Demongeot, B.~Lacolle (eds.), NumericaL Methods in the Study of
  Critical Phenomena, Springer-Verlag, Berlin, 1981, pp. 12 -- 25.

\bibitem{Simon/1970a} B.~Simon, Coupling constant analyticity for the
  anharmonic oscillator, Ann.  Phys. (NY) 58 (1970) 76 -- 136.

\bibitem{Simon/1991} B.~Simon, Fifty years of eigenvalue perturbation
  theory, Bull. Amer. Math. Soc.  24 (1991) 303 -- 319.

\bibitem{Smith/Ford/1979} D.~A. Smith, W.~F. Ford, Acceleration of linear
  and logarithmic convergence, SIAM J. Numer. Anal. 16 (1979) 223 -- 240.

\bibitem{Smith/Ford/1982} D.~A. Smith, W.~F. Ford, Numerical comparisons
  of nonlinear convergence accelerators, Math. Comput. 38 (1982) 481 --
  499.

\bibitem{Srivastava/Choi/2001} H.~M. Srivastava, J.~Choi, Series
  Associated with the Zeta and Related Functions, Kluwer, Dordrecht,
  2001.

\bibitem{Thomann/1990} J.~Thomann, Resommation des series formelles,
  Numer. Math. 58 (1990) 503 -- 535.

\bibitem{Thomann/1995} J.~Thomann, Proc\'{e}d\'{e}s formels et
  num\'{e}riques de sommation de s\'{e}ries d'\'{e}quations
  differentielles, Expo. Math. 13 (1995) 223 -- 246.

\bibitem{Tweddle/2003} I.~Tweddle, James Stirling's Methodus
  Differentialis: An Annotated Translation of Stirling's Text,
  Springer-Verlag, London, 2003.

\bibitem{Watson/1912b} G.~Watson, The transformation of an asymptotic
  series into a convergent series of inverse factorials, Rend. Circ. Mat.
  Palermo 34 (1912) 41 -- 88.

\bibitem{Weniger/1989} E.~J. Weniger, Nonlinear sequence transformations
  for the acceleration of convergence and the summation of divergent
  series, Comput. Phys. Rep. 10 (1989) 189 -- 371, {L}os Alamos Preprint
  math-ph/0306302 (\texttt{http://arXiv.org}).

\bibitem{Weniger/1990} E.~J. Weniger, On the summation of some divergent
  hypergeometric series and related perturbation expansions, J. Comput.
  Appl. Math. 32 (1990) 291 -- 300.

\bibitem{Weniger/1992} E.~J. Weniger, Interpolation between sequence
  transformations, Numer. Algor. 3 (1992) 477 -- 486.

\bibitem{Weniger/1994b} E.~J. Weniger, Verallgemeinerte
  {S}ummationsprozesse als numerische {H}ilfsmittel f\"{u}r
  quantenmechani\-sche und quantenchemische {R}echnungen, {H}abilitation
  thesis, Fachbereich Chemie und Pharmazie, Universit\"at Regensburg,
  {L}os Alamos Preprint math-ph/0306048 (\texttt{http://arXiv.org})
  (1994).

\bibitem{Weniger/1996c} E.~J. Weniger, {A} convergent renormalized strong
  coupling perturbation expansion for the ground state energy of the
  quartic, sextic, and octic anharmonic oscillator, Ann. Phys. (NY) 246
  (1996) 133 -- 165.

\bibitem{Weniger/1996e} E.~J. Weniger, {C}onstruction of the strong
  coupling expansion for the ground state energy of the quartic, sextic
  and octic anharmonic oscillator via a renormalized strong coupling
  expansion, Phys. Rev. Lett. 77 (1996) 2859 -- 2862.

\bibitem{Weniger/2004} E.~J. Weniger, Mathematical properties of a new
  {L}evin-type sequence transformation introduced by
  {{\v{C}}\'{\i}\v{z}ek}, {Z}amastil, and {Sk\'{a}la}. {I}. {A}lgebraic
  theory, J. Math. Phys. 45 (2004) 1209 -- 1246.

\bibitem{Weniger/2007a} E.~J. Weniger, Asymptotic approximations to
  truncation errors of series representations for special functions, in:
  A.~Iske, J.~Levesley (eds.), Algorithms for Approximation,
  Springer-Verlag, Berlin, 2007.

\bibitem{Weniger/2008} E.~J. Weniger, On the analyticity of {L}aguerre
  series, J. Phys. A 41 (2008) 425207--1 -- 425207--43.

\bibitem{Weniger/Cizek/Vinette/1993} E.~J. Weniger,
  J.~{\v{C}\'{\i}\v{z}ek}, F.~Vinette, The summation of the ordinary and
  renormalized perturbation series for the ground state energy of the
  quartic, sextic, and octic anharmonic oscillators using nonlinear
  sequence transformations, J. Math. Phys. 34 (1993) 571 -- 609.

\bibitem{Weniger/Steinborn/1989a} E.~J. Weniger, E.~O. Steinborn,
  Nonlinear sequence transformations for the efficient evaluation of
  auxiliary functions for {GTO} molecular integrals, in: M.~Defranceschi,
  J.~Delhalle (eds.), Numerical Determination of the Electronic Structure
  of Atoms, Diatomic and Polyatomic Molecules, Kluwer, Dordrecht, 1989,
  pp. 341 -- 346.

\end{thebibliography}

\end{document}